\documentclass[reqno]{amsart}
\usepackage{hyperref} %Package hyperref+amsart : use to make link for index

\begin{document}
\title[\hfilneg \hfil The non-uniqueness of solutions of  Navier-Stokes equations.]
{Non-uniqueness of solutions for 3D  Navier-Stokes equations in bounded domains.}

\author[Vu Thanh Nguyen \hfil \hfilneg]
{Vu Thanh Nguyen}

\address{Vu Thanh Nguyen 
	\newline
	Department of Mathematics, University of
	Texas at San
	Antonio.
\newline  \hspace*{.2 in}{\it  Family name}: $NGUY\widetilde{\hat{E}}N$, {\it Given name}: $V\widetilde{U}$.
	}
\email{vu.nguyen@utsa.edu,~~or 9vunguyen@gmail.com}

\thanks{  }

\begin{abstract}
This paper examines the uniqueness/non-uniqueness of local-in-time strong solutions for the incompressible 3D Navier-Stokes equations in bounded domains, which  are  $\partial_t u=\nu \Delta u- u\cdot \nabla u-\nabla p+ f$ and $div~u=0$.

The focus of this study is on the case where the boundary condition is defined as $u\cdot \vec{n}|_{\partial\Omega}=0$, indicating that the fluid is unable to pass through the boundary but is allowed to move tangentially to it. By presenting a counterexample, this paper demonstrates the existence of two distinct strong solutions to the Navier-Stokes equations under this boundary condition.

\end{abstract}

\maketitle	
	%%%%%%%%%%%%%%%%%%
		\newcommand{\la}{\langle}
		\newcommand{\ra}{\rangle}
\newcommand{\dps}{\displaystyle}
\newcommand{\pg}{\partial}
\newcommand{\ino}{ \int_\Omega} 
\newcommand{\RR}{\mathbb{R}}
\newcommand{\NN}{\mathbb{N}}
\newcommand{\LL}{\mathbb{L}}
\newcommand{\HH}{\mathbb{H}}
\newcommand{\hh}{\hspace*{.2in}}
\numberwithin{equation}{section} % phuong trinh se duoc ghi theo tung section, khong lien tuc ca bai  
\numberwithin{equation}{subsection}
%\allowdisplaybreaks
\newtheorem{theorem}{Theorem}[subsection]
\newtheorem{defi}[theorem]{Definition}
\newtheorem{lemma}[theorem]{Lemma}
\newtheorem{rem}[theorem]{Remark}
\newtheorem{exam}[theorem]{Example}
\newtheorem{propo}[theorem]{Proposition}
\newtheorem{corol}[theorem]{Corollary}
\newtheorem{intro}[theorem]{Introduction}
\newtheorem{Notations}[theorem]{Notations}
\newtheorem{claim}[theorem]{Claim}
%\tableofcontents
\allowdisplaybreaks
\section{Introduction}  
 \hspace*{-0.2in}$-$  
 	We consider the incompressible Navier-Stokes equations: 
 	%x----------3101
 	\begin{equation}\label{3101}
 	(\mbox{NS})\left\{
 	\begin{array}{rll}
 	\partial_t u&=\nu \Delta  u- u\!\cdot\!  \nabla  u-\nabla p+ f , &\\
 	div\,u&=0,&\\
 	u(.,0)& = u_o,&
 	\end{array}
 	\right. 
 	\end{equation}
 in $\Omega\times [0,+\infty)$, where $\Omega$ is a bounded domain in $\RR^3$, the force $f$ and the initial velocity $u_o$ are  given.
 	
 %%%%%%%%%%
%%%%%%%%%%%%%%%
If $u$ solves the problem (NS), then $u$ satisfies

\begin{equation}\label{3134}
\int_{\partial\Omega}u\cdot \vec{n}\, d\sigma = \int_\Omega div\,u \,dx = \int_\Omega \!0 \,dx = 0
\end{equation}
by the Gauss divergence theorem, where $\vec{n}$ is the normal vector to the boundary $\partial\Omega$. Hence,  boundary conditions for the problem (NS) must satisfy (\ref{3134}). We observe that the natural boundary condition satisfying (\ref{3134}) is
\begin{align}\label{3102}
	u\cdot \vec{n}=0\text{ on }\partial\Omega,
\end{align}
which  indicates  that the fluid cannot pass through the boundary but is allowed to move tangentially to it.

The uniqueness of local-in-time strong solutions to the problem (NS) under the boundary condition $u|_{\partial\Omega}=0$, which is a special case of the condition described by equation (\ref{3102}), has been extensively studied in various books such as \cite{Lady}, \cite{Constantin}, \cite{Boyer}, \cite{Robinson}, \cite{Teman}, and others. However, the uniqueness/non-uniqueness of local-in-time  strong solutions to the problem (NS) under the boundary condition (\ref{3102}) has not yet been established. 

In other words, the uniqueness of local-in-time  strong solutions  in the space $C([0, \epsilon];H\cap (H_o^1(\Omega))^3)$ to the problem (NS) has been studied in various books, where
$ \dps H:=\{ v\in (L^2(\Omega))^3: div\,v=0 \mbox{ in }\Omega, \mbox{ and } v\cdot \vec{n}|_{\pg\Omega}=0 \}$.
However, the uniqueness/non-uniqueness of local-in-time  strong solutions  in the space $C([0, \epsilon];H\cap (H^1(\Omega))^3)$ to  the problem (NS) has not yet been established.
This paper aims to investigate this uniqueness/non-uniqueness.

The main result of this paper, which will be proved in the last section (Section 3), is stated as follows:

\begin{theorem}\label{3128}
	Let $\Omega$ be a ball centered at the origin in $\mathbb{R}^3$. Suppose $T>0$, $u_o\in V\cap (H^2(\Omega))^3\cap C(\bar{\Omega})$, and $f\in L^2(0,T;(L^2(\Omega))^3)$.
	
	 Then there exists $\epsilon\in (0, T]$ such that the problem (NS) with the boundary condition (\ref{3102}) has at least two distinct strong solutions on the time interval $[0, \epsilon]$.
\end{theorem}

The result in Theorem \ref{3128} demonstrates that the uniqueness of local-in-time strong solutions to the problem (NS) with the boundary condition (\ref{3102}) is not valid.

In the next section, we will provide the definition of strong solutions for solutions in Theorem \ref{3128}.  Additionally, we will present a result that pertains to the existence of local-in-time strong solutions to a problem, which is  similar to Navier-Stokes equations  under  the boundary condition $u|_{\partial\Omega}=0$.
%%%%%%%%%%%%%%%%%
\section{Preliminary results}\label{3027}
\begin{Notations}\label{3123}~ The following notations will be used through this paper:
% Let $\Omega\subset \RR^3$ be a bounded domain of class $C^2$.
\begin{itemize}
%\item[\tiny$\bullet$ ]  	Let $\Omega$ be a bounded domain $\Omega$ in $\RR^3$ of class $C^3$.
\item[\tiny$\bullet$ ] $\LL^2:=(L^2(\Omega))^3$, $\HH^s:=(W^{s,2}(\Omega))^3$, where $W^{s,2}(\Omega)$ is a Sobolev space.
\item [\tiny$\bullet$ ] $\|.\|_q:=\|.\|_{\LL^q(\Omega)}$.
\item[\tiny$\bullet$ ] $u\cdot\nabla v:=(u\cdot \nabla)v$.
 \item[\tiny$\bullet$ ] $ \dps H:=\{ v\in \LL^2: div\,v=0 \mbox{ in }\Omega, \mbox{ and } v\cdot \vec{n}\,|_{\pg\Omega}=0 \},$ where $\vec{n}$ is the normal vector to the boundary $\pg\Omega$.
\item[\tiny$\bullet$ ] $ \dps V:=\{ v\in \HH_0^1: div\, v=0\}=H\cap \HH_0^1$.
%   where $\HH_o^1(\RR^3)=\HH^1(\RR^3)$.        
% $-$ $V:=\{u\in (H_o^1(\Omega))^3: div\,u=0\}$. \\
\item[\tiny$\bullet$ ]  The space of test functions on the space-time domain $\Omega\times [0,T]$ is given by
%x-----3154,55
\begin{align}\label{3154}
&\mathcal{D}_{\sigma,T}=\{\varphi\in C_c^\infty(\Omega\times [0, T)]: div~\varphi(t)=0 \mbox{ for all }t\in [0,T]\}.
%\\
%\label{3155}
%& \tilde{\mathcal{D}}_{\sigma,T}=\{\varphi: \varphi=\sum_{k=1}^N\alpha_k(t)a_k(x), \alpha_k\in C_c^1([0,T]), %a_k\in\mathcal{N}, N\in\NN\},
\end{align}
%where $\mathcal{N}$  represents the basis of $H$ consisting of eigenfunctions of the Stokes operator. For more detailed %information, please refer to  Robinson \& Rodrigo\& Sadowski \cite{ Robinson \& Rodrigo\& Sadowski}.
%\item[\tiny$\bullet$ ] The curl of a vector filed $(g_1, g_2, g_3)$ is defined by
%%x-----3130
%\begin{equation}\label{3130}
%\nabla\times (g_1, g_2, g_3):=\left(\frac{\pg g_3}{\pg x_2}-\frac{\pg g_2}{\pg x_3}\,,\,
%                \frac{\pg g_1}{\pg x_3}-\frac{\pg g_3}{\pg x_1}\,,\,
%                \frac{\pg g_2}{\pg x_1}-\frac{\pg g_1}{\pg x_2}\right)
%\end{equation}  
%\item[\tiny$\bullet$ ] $|x|:=\sqrt{x_1^2+x_2^2+x_3^2}$
\end{itemize}
%%%%%%%%%%%%%%%%
\end{Notations}%%%%%

%%%%%%%%%%%%%%%%

Under the boundary condition $u|_{\pg\Omega}=0$, strong solutions belong to  $C([0,T];V)$, where $V=H\cap \HH_o^1$ (prefer to Definition of strong solutions in  Robinson \& Rodrigo \& Sadowski \cite{Robinson}).
However, in the definition below, the boundary condition is  $u\cdot\vec{n}\,|_{\pg\Omega}=0$ instead of  $u|_{\pg\Omega}=0$, hence  the  space for solutions will be $C([0,T], H\cap \HH^1)$ instead of $C([0,T];V)$.
%x----3127
\begin{defi}\label{3127} Let $\Omega\subset\RR^3$ be a bounded domain with smooth boundary.
	
	We say that a function $u$ is a strong solution  of the problem (NS) with the boundary condition   $u\cdot\vec{n}\,|_{\pg\Omega}=0$
	%with the boundary condition (\ref{3102})
	on the time interval $[0, T]$ if $u $ satisfies the following:
	%x---3173, 3174,3114
	\begin{align}\label{3173}
	&\mbox{\tiny$\bullet$ } u|_{t=0}=u_o,\\ \label{3174}
	&\mbox{\tiny$\bullet$ } u\in C([0,T];H\cap \HH^1)\cap L^2(0, T; \HH^2), ~ \pg_tu\in L^2(0,T;\LL^2), \\
	\label{3114}
			&\mbox{\tiny$\bullet$ } 
	  	\int_0^s\la \pg_tu-\Delta u+(u\cdot\nabla)u-f,\,\varphi\ra\, dt=0,
	\end{align}	
	\hh	\hh  	for all test functions $\varphi\in \mathcal{D}_{\sigma,T}$ and almost every $s\in [0,T]$.

\end{defi}
\begin{rem}~
	\\
	$-$ Due to Robinson \cite[Proposition 1.35]{Robinson}, we can deduce  from the two facts that $u\in  L^2(0, T; \HH^2) $ and $\dps \pg_tu\in L^2(0,T;\LL^2)$ that $u\in C([0, T];  \HH^1)$. Hence, the property $u\in C([0, T];  \HH^1)$ can be omitted from Definition \ref{3127}.\\
	$-$  If $u$ satisfies (\ref{3174}), then we can deduce from the properties of the space $H$ that  $u\cdot\vec{n}\,|_{\pg\Omega}=0$, and  $div\,u=0$ in $\Omega$. 
\end{rem} 
%%%%%%%%%%%%%%%%%%%%%
The following inequality will   be employed in the proof of next theorem (Theorem \ref{3110}).
\begin{lemma}\label{5140}(Gagliardo-Nirenberg's inequality)
	
	Let $\Omega$  be a ball in $\RR^3$ and  $ q\in [2,6]$. 
	There exists a  constant   $C_1, C_2$  such that
	\begin{align*}
		&\mbox{\tiny $\bullet$ }\dps	\|u\|_q \le C_1 \|u\|_2^{3/q-1/2}\|u\|_{\HH^1}^{3/2-3/q}
		,~\forall  u\in  \HH^{1},
		\\
		& \mbox{\tiny $\bullet$ }\dps \|u\|_q \le C_2 \|u\|_2^{3/q-1/2}\|\nabla u\|_2^{3/2-3/q},~ \forall  u\in  \HH_o^1.
	\end{align*}	 
\end{lemma}~
{\bf Proof.}  

The inequalities  can be derived  from  Boyer \& Fabrie \cite[Proposition III.2.35, Remark III.2.17]{Boyer}. 
%$-$ Based on inequality (\ref{3141}) and the fact that $\|u\|_{\HH^2}\leq c\|Au\|_2$, we can deduce the inequality (\ref{3143}).
\hfill $\Box$

The system of equations in the following statement is similar to the system of Navier-Stokes equations. The existence of a strong solution for this problem is analogous to the existence of a strong solution for the Navier-Stokes equations under the boundary condition $u|_{\pg\Omega}=0$.    
%x-----5140

%x---3110
\begin{theorem}\label{3110} 	Let $\Omega\subset\RR^3$ be a ball and  $T>0$. Suppose  $u_o\in V\cap \HH^2$, $\beta \in C([0,T];\HH^2)$ and $F\in L^2(0,T;\LL^2 )$. 
	
	Consider the equations for unknown $(u,p)$ in $\Omega\times [0,T]$: 
	%x----------3007
	\begin{equation}\label{3007}
	\left\{
	\begin{array}{rll}
	\partial_t u&=\nu \Delta  u- u\!\cdot\!  \nabla  u- \beta \!\cdot\!  \nabla  u- u\!\cdot\!  \nabla  \beta -\nabla p+ F, & \\
		div\,u&=0,&\\
			u(.,0)& = u_o.&\\
	\end{array}
	\right. 
	\end{equation}
	
  The boundary  condition is 
  %x-----3106
  \begin{align}\label{3106}
  	u|_{\pg\Omega}=0.
  \end{align}

  Then, there exists $T_1\in (0, T]$ such that   the problem $(\ref{3007})$ possesses a  solution $u$ that satisfies: 
  %x----- 3172,3166, 3167
  \begin{align}
  	 \label{3172}
  	 &\mbox{\tiny $\bullet$ } u|_{t=0}=u_o,\\
  	\label{3166}
  	&\mbox{\tiny $\bullet$ } u\in C([0,T];V)\cap L^2(0, T; \HH^2), ~ \pg_tu\in L^2(0,T;\LL^2),\\
  	\label{3167}
  	&\mbox{\tiny $\bullet$ } 
  \int_0^s\la 	\partial_t u-\nu \Delta  u+ u\!\cdot\!  \nabla  u+ \beta \!\cdot\!  \nabla  u+ u\!\cdot\!  \nabla  \beta - F,\,\varphi\ra\, dt=0  	
  	\end{align}
 \hspace{.8in}for all test functions $\varphi\in \mathcal{D}_{\sigma,T_1}$ and almost every $s\in [0,T_1]$.
\end{theorem}
\begin{rem}~\\
$-$ In the special case where $\beta\equiv 0$, problem (\ref{3007}) simplifies to the Navier-Stokes equations with the boundary condition (\ref{3106}), which has been extensively studied in various books. 	
\\
$-$
Both additional terms $(\beta\cdot\nabla)u$ and $(u\cdot\nabla)\beta$ are only linear with respect to the unknown $u$. Therefore, the proof of Theorem \ref{3110} follows a similar approach to that used for the Navier-Stokes equations with the same  boundary condition.
\end{rem}
%%%%%%%%%%
%%%%%%%%%%
  {\bf Proof}.

This theorem will be proven using the Galerkin method, which follows a similar approach to that used for the Navier-Stokes equations. The sketch of the proof is as follows:
\begin{enumerate}
\item[\mbox{\tiny $\bullet$} \it Step 1:] {\it    Galerkin approximations}.

 The {\it n}-th order Galerkin problem corresponding to the problem (\ref{3007}) is given by
%x-----3135
\begin{align}\label{3135}
	\hspace*{.4in}		\left\lbrace 
	\begin{array}{rl}
		\displaystyle \partial_t u _n+A u_n&=\mathcal{P}_n(- (u _n\!\cdot\!\nabla) u _n -   (\beta\!\cdot\!\nabla) u _n - (u _n\!\cdot\!\nabla)     \beta  +   F   ), \\
		u _n(.,0)&=\mathcal{P}_n( u _o), \\
		div\,  u _n&=0,	
	\end{array}
	\right.
\end{align}
where $A$ is the Stokes operator. Here,
$\mathcal{P}_n$ is the projection operator $\LL^2\to H$ defined by 
$\dps P_nu=\sum_{i=1}^n\la u,a_i\ra a_i$,  $\dps u_n=\sum_{i=1}^n \alpha_k(t)a_k$,
%$ \la u_n, a_i\ra$,
%orthonormal basis in $H$ and an orthogonal basis in $V$, 
 where $\{a_1, a_2,...\}$ is the family of eigenfunctions of the Stokes operator. 

Taking the $L^2$-inner product the first equation in (\ref{3135})  with $a_k$ for $1\le k\le n$, we  get    $n$ ordinary differential equations with $n$ unknowns $(\alpha_k(t))_{1\le i\le n}$. This   ODE system is quadratic with  coefficients in $C^1[0,T]$. If we denote the vector of components $(\alpha_k(t))_{k\le n}$ as $\alpha(t)$, the ODE system can be written as the form
\begin{equation}
\frac{d\alpha}{dt}=\mathcal{F}(t,\alpha),
\end{equation}
where $\mathcal{F}$ is continuous and locally Lipschitz continuous with respect to $\alpha$. 

  The Cauchy-Lipschitz theorem gives us the existence of a unique $C^1$ solution, which is defined on a maximum interval belonging to $[0,T]$. Besides, by the estimate of $\|u_n\|_2^2$ at (\ref{3136}) in the next step, 
  $\dps \sum_{i=1}^n|\alpha_i(t)|^2$ is defined on $[0, T]$. So that $\alpha_k\in C^1([0,T])$ for $1\le k\le n$.
Hence,  the problem (\ref{3135}) has a unique solution  $u_n\in C^1([0, T]; V\cap \HH^2)$. 
%%%%%%%%%%%%%%%%%%%%%%%
\item[\mbox{\tiny $\bullet$} \it Step 2:]  {\it The boundedness of sequences  $(\|u_n\|_{L^\infty(0,T;\LL^2)}^2)_n$, $(\int_0^{T}\!\!\|\nabla u _n\|_2^2\, dt)_n$. }
\\
$-$ We take the $L^2$-inner product of the first equation in (\ref{3135}) with $u_n$, and then utilize H\"older's inequality  and Gagliardo-Nirenberg's inequality. This yields the following inequality for every $t\in [0, T]$ and $n\in\mathbb{N}$: 
\begin{align*}
	\hh &\frac{d}{dt}\| u _n\|_2^2+2\nu \|\nabla  u _n\|_2^2\le \\
	&\le  0+2\|\beta\|_6\|\nabla u_n\|_2\|u_n\|_3+2\|u_n\|_4\|\nabla \beta\|_2\|u_n\|_4+2\|  F   \|_2\| u _n\|_2\\
	&\le c_1\|  \beta  \|_{H^1}\|\nabla u_n\|_2\| u _n\|_2^\frac{1}{2}\|\nabla u _n\|_2^\frac{1}{2} 
	+ c_1\| \nabla    \beta  \|_2\| u _n\|_2^\frac{1}{2}\|\nabla u _n\|_2^\frac{3}{2} 
	+2\|  F   \|_2\| u _n\|_2\\
	&\le c_2\| u _n\|_2^\frac{1}{2}\|\nabla u _n\|_2^\frac{3}{2} 
	+c_2\| u _n\|_2.
\end{align*} 

Applying Young's inequality, we obtain:
\begin{align*}
	\frac{d}{dt}\| u _n\|_2^2+2 \nu \|\nabla  u _n\|_2^2
	&\le
	(c_3\| u _n\|_2^2+\nu\|\nabla u _n\|_2^2)+ (\frac{1}{2}\|u _n\|_2^2+\frac{1}{2}c_2^2).
\end{align*}

Therefore, 
%x-------------3137
\begin{equation}\label{3137}
\frac{d}{dt}\| u _n\|_2^2+ \nu \|\nabla  u _n\|_2^2\le c_4\| u _n\|_2^2+c_4.
\end{equation}

Omitting  the second term, we have a  Gronwall's inequality. So that
\begin{align*}
	\forall t\in [0, T], n\in\mathbb{N}, \| u _n\|_2^2& \le e^{c_4T}(\| u _n(.,0)\|_2^2+c_4T)\\
	& \le e^{c_4T}(\|u _o\|_2^2+c_4T)=:m_o.
\end{align*}

Hence, the sequence $(\|u_n\|_2^2)_n$ is bounded uniformly on $[0, T]$:
%x----3136
\begin{align}\label{3136}
	\forall t\in [0, T],\forall  n\in\mathbb{N}, \| u _n\|_2^2& \le m_o.
\end{align}
$-$ 
Integrating both sides of (\ref{3137})  in $t$ between $0$ and $T$, and using (\ref{3136}), we deduce that 
\begin{align*}
&	\forall n\in\mathbb{N}, \| u _n(.,T)\|_2^2-\|u_n(.,0)\|_2^2\!+\nu\int_0^{T}\!\!\|\nabla u _n\|_2^2\, dt \!\le 
	(c_4m_0+c_4)T\\
&	\forall n\in\mathbb{N},~ \nu\int_0^{T}\!\!\|\nabla u _n\|_2^2\, dt \!\le \|u_n(.,0)\|_2^2+c_5\le 
\|u_o\|_2^2+c_5 .
\end{align*}

Hence, the sequence $(\int_0^{T}\!\!\|\nabla u _n\|_2^2\, dt)_n$ is bounded:
%x------3138
\begin{align}\label{3138}
	\forall n\in\mathbb{N},\int_0^{T}\!\!\|\nabla u _n\|_2^2\, dt \!&\le	c_6.
\end{align}

%%%%%%%%%%%%%%%%%%%%%%%
\item[\mbox{\tiny $\bullet$} \it Step 3:]  {\it The uniform boundedness of the sequence $(\|\nabla u_n\|_2^2)_n$ }.

%This step in important, because from the uniform boundedness of    $(\|\nabla u_n\|_2)_n$ we can implies the uniform %boundedness of other sequences in next steps.
 
Taking the $L^2$-inner product  the   first equation in  (\ref{3135}) by $A u _n$ and using H\"older's inequality,  we obtain: 
\begin{align*}
	&\frac{d}{dt} \|\nabla u_n\|_2^2+2\nu\|A u_n\|_2^2
		\le
	2\| u_n\|_6\|\nabla u_n\|_3\|A u_n\|_2
	+
	\\&\hspace*{.3in}+2\|   \beta  \|_6\|\nabla u_n\|_3\|A u_n\|_2+2\| u_n\|_3\|\nabla   \beta  \|_6\|A u_n\|_2 +
	  2\|  F   \|_2\|A u_n\|_2.
	% % % % % % % % % % % % % % % % % %  3
\end{align*}

We use Gagliardo-Nirenberg's inequality, along with  (\ref{3136}) and the fact 	$\|u_n\|_{\HH^2(\Omega)}\le c\|Au_n\|_{\LL^2(\Omega)}$ to derive the following:
\begin{align*}
	& \frac{d}{dt} \|\nabla u_n\|_2^2+2\nu\|A u_n\|_2^2
	\le c_7
	\|\nabla u_n\|_2\|\nabla u_n\|_2^\frac{1}{2}\|A u_n\|_2^\frac{1}{2}\!
	\|A u_n\|_2+\\
	&\hspace*{0.1in}+ c_7\|\beta\|_{\HH^1} \|\nabla u_n\|_2^\frac{1}{2}\|A u_n\|_2^\frac{1}{2}\|A u_n\|_2
	\!+\!   c_7\| u_n\|_2^\frac{1}{2}\!  \|\nabla v_n\|_2^\frac{1}{2}\|\beta\|_{\HH^2}\|A u_n\|_2+\\
	&\hspace*{1in}+  c_7\|A u_n\|_2
	\end{align*}
	\begin{align*}	
	& \frac{d}{dt} \|\nabla u_n\|_2^2+2\nu\|A u_n\|_2^2
	\le\,
	c_7\|\nabla u_n\|_2^\frac{3}{2}\|A u_n\|_2^\frac{3}{2}\!
	+\\
	&\hspace*{0.6in}+  c_8\|\nabla u_n\|_2^\frac{1}{2}\|A u_n\|_2^\frac{3}{2}
	\!+\!  c_8 \|\nabla u_n\|_2^\frac{1}{2}\|A u_n\|_2
	\!+c_7\|A u_n\|_2.
\end{align*}

Using Young's inequality yields 
\begin{align*}	
	& \frac{d}{dt} \|\nabla u_n\|_2^2+2\nu\|A u_n\|_2^2
	\le (c_9
	\|\nabla u_n\|_2^6+ \frac{\nu}{4}\|A u_n\|_2^2)+\\
	&\hspace*{.5in}+  (c_9\|\nabla u_n\|_2^2+ \frac{\nu}{4}\|A u_n\|_2^2)
	\!+\!  ( c_9 \|\nabla u_n\|_2+ \frac{\nu}{4}\|A u_n\|_2^2)+\\
	&\hspace*{1in}\!+(c_9+ \frac{\nu}{4}\|A u_n\|_2^2)\\
		& \frac{d}{dt} \|\nabla u_n\|_2^2+\nu\|A u_n\|_2^2
		\le c_9
		\|\nabla u_n\|_2^6+  c_9\|\nabla u_n\|_2^2
		\!+\!   c_9 \|\nabla u_n\|_2+ c_9.
\end{align*}

Here, $\|\nabla u_n\|_2^6\le (\|\nabla u_n\|^2+1)^3$,
 $\|\nabla u_n\|_2^2\le (\|\nabla u_n\|_2^2+1)^3$, and \\
 $\|\nabla u_n\|_2\le (\|\nabla u_n\|^2+1)^3$. 
 Therefore, 
%x-----3422
\begin{align}\label{3422}
\forall n\in\NN,	\frac{d}{dt} \|\nabla u_n\|_2^2+\nu \|A u_n\|_2^2
	\le c_{10}(\|\nabla u_n\|_2^2+1)^3.   
\end{align}

By omitting the second term and denoting  $X_n:=\|\nabla u_n\|_2^2+1$, we get
$\dps X_n(0)=( \|\nabla u_n(.,0)\|_2^2+1)\le \|u_o\|_{\HH^1}^2+1$ and
\begin{align*}
\forall n\in\NN,\hh	\frac{dX_n}{dt}  &\le c_{10}X_n^3
\\
	\frac{d X_n}{dt}X_n^{-3}& \le c_{10}.
\end{align*}

Integrating both sides in time between $0$ and $t$, we obtain 
\begin{align*}
\forall n\in\NN, ~X_n(0)^{-2}- X_n(t)^{-2}&\le 2c_{10}t\\
%%	X(0)- \left( \|\nabla u_n(.,t)\|_2^2+1\right) ^{-2}	&\le  2c_9 t\\
	~X_n(t)^{-2}	&\ge  X_n(0)^{-2}- 2c_{10} t\ge (\|u_o\|_{\HH^1}^2+1)^{-2} - 2c_{10} t\\
~X_n(t)^{-2}&\ge\frac{1}{2} (\|u_o\|_{\HH^1}^2+1)^{-2}	\mbox{ for }t\le \frac{(\|u_o\|_{\HH^1}^2+1)^{-2}}{4c_{10}}.
\end{align*}

Choosing $\dps 
T_1:=\min\left\{\frac{(\|u_o\|_{\HH^1}^2+1)^{-2}}{4c_{10}}; T\right\},
$
we get 
\begin{align*}
	\forall t\in [0,T_1], \forall n\in\NN, 
	X_n(t)^{-2}&\ge \frac{1}{2} (\|u_o\|_{\HH^1}^2+1)^{-2} \\
	~\left( \|\nabla u_n(.,t)\|_2^2+1\right) ^{-2}	&\ge \frac{1}{2} (\|u_o\|_{\HH^1}^2+1)^{-2} \\
	~ \|\nabla u_n(.,t)\|_2^2+1	&\le \sqrt{2}( \|u_o\|_{\HH^1}^2+1).
\end{align*}

From  this fact and (\ref{3136}), we can conclude that the sequence \\
$(\| u_n(.,t)\|_{\HH^1}^2)_n$ is bounded uniformly on $[0, T_1]$:
%x---------------3441
\begin{equation}\label{3441}
\exists m_1>0, \forall t\in [0,T_1], \forall n\in\NN, ~ \| u_n(.,t)\|_{\HH^1}^2	\le m_1. 
\end{equation}
\item[\mbox{\tiny $\bullet$} \it Step 4:]  {\it The  boundedness of the sequence $ (\int_0^{T_1}\! \| u _n\|_{\HH^2}^2 \,dt)_n$. }

By integrating both sides of equation (\ref{3422}) in time  between $0$ and $T_1$, and using (\ref{3441}), we obtain the following result for every $n\in\NN$:
\begin{align*}
&	 \|\nabla u_n(.,T_1)\|_2^2-\|\nabla u_n(.,0)\|_2^2+\nu\int_0^{T_1} \|A u_n\|_2^2\,dt
		\le c_{10}(m_1+1)^3T_1   \\
&\nu 	\int_0^{T_1} \|A u_n\|_2^2\,dt
		\le \|\nabla u_n(.,0)\|_2^2+c_{10}(m_1+1)^3T_1\le m_1+c_{10}(m_1+1)^3T_1 .   
\end{align*}	

Therefore, the sequence $(\int_0^{T_1} \|A u_n\|_2^2\,dt)_n$ is bounded. From this boundedness and the fact $\|u_n\|_{\HH^2}\le  c\|Au_n\|_2$, we deduce that the sequence $(\int_0^{T_1} \| u_n\|_{\HH^2}^2\,dt)_n$ is bounded:
%x----3139
\begin{equation}
\label{3139}
\exists m_2>0,~	\forall n\in\NN,~	\int_0^{T_1}\| u_n\|_{\HH^2}^2\,dt
	\le m_2
\end{equation}
\item[\mbox{\tiny $\bullet$} \it Step 5:]  {\it  The  boundedness of the sequence $ ( \int_0^{T_1}\|(u_n\cdot\nabla)u_n\|_2^2\,dt)_n$. }

By using  the  boundedness of the sequences $(\|u_n\|_{L^\infty(0,T_1;\HH^1)})_n$ and\\
  $ (\int_0^{T_1}\! \| u _n\|_{\HH^2}^2 \,dt)_n$ (as stated at (\ref{3441}) and (\ref{3139})),  we  get
\begin{align*}
	 \forall n\in\NN,~ \int_0^{T_1}	\|(u_n\cdot\nabla)u_n\|_2^2\,dt &\le  \int_0^{T_1}\ino |u_n|^2|\nabla u_n|^2 dx\,dt\\
	&\le \int_0^{T_1}\ino |\nabla u_n|^2 dx\|u_n(.,t)\|_{\LL^\infty}^2\,dt \\
		&\le \|u_n\|^2_{L^\infty(0,T_1;\HH^1)} \int_0^{T_1}\|u_n(.,t)\|_{\LL^\infty}^2\,dt \\
	&\le c_{11}\|u_n\|^2_{L^\infty(0,T_1;\HH^1)} \int_0^{T_1} \|u_n(.,t)\|_{\HH^2}^2\,dt\\
	&\le c_{11}m_1m_2,
\end{align*}	
where we have used the embedding of $\HH^2(\Omega)$ into $\LL^\infty(\Omega)$.
\item[\mbox{\tiny $\bullet$} \it Step 6:]  {\it The  boundedness of the sequence $ (\int_0^{T_1}\|(\beta\cdot\nabla)u_n\|_2^2\,dt)_n$ }.

Similarly to Step 5, based on  the  boundedness of $(\|u_n\|_{L^\infty(0,T_1;\HH^1)})_n$
 (as stated at (\ref{3441})), and  the assumption $\beta\in C([0, T];\HH^2)$,   we  obtain
\begin{align*}
\forall n\in\NN,~ \int_0^{T_1}	\|(\beta\cdot\nabla)u_n\|_2^2\,dt &\le  \int_0^{T_1}\ino |\beta|^2|\nabla u|^2 dx\,dt\\
&\le \int_0^{T_1}\ino |\nabla u|^2 dx\|\beta(.,t)\|_{\LL^\infty}^2\,dt \\
&\le\|u\|^2_{L^\infty(0,T_1;\HH^1)} \int_0^{T_1} \|\beta(.,t)\|_{\LL^\infty}^2\,dt\\
&\le c_{12}\|u\|^2_{L^\infty(0,T_1;\HH^1)} \int_0^{T_1} \|\beta\|_{\HH^2}^2\,dt\\
&\le c_{12} m_1 \int_0^{T_1} \|\beta\|_{\HH^2}^2\,dt.
\end{align*}	

\item[\mbox{\tiny $\bullet$} \it Step 7:]  {\it The boundedness of the sequence $ (\int_0^T\|(u_n\cdot\nabla)\beta\|_2^2\,dt)_n$.}

Similarly to Step 5, based on the  boundedness of  $(\int_0^{T_1}\|u_n\|_{\HH^2}^2)_ndt$ 
(as stated at (\ref{3139}), and  the assumption $\beta\in C([0, T];\HH^2)$,   we obtain
\begin{align*}
	\forall n\in\NN,~ \int_0^{T_1}	\|(u_n\cdot\nabla)\beta\|_2^2\,dt &\le  \int_0^{T_1}\ino |u_n|^2|\nabla \beta|^2 dx\,dt\\
	&\le \int_0^{T_1}\ino |\nabla \beta|^2 dx\|u_n(.,t)\|_{\LL^\infty}^2\,dt \\
	&\le c_{13}\|\beta \|^2_{L^\infty(0,T_1;\HH^1)} \int_0^{T_1} \|u_n\|_{\HH^2}^2\,dt\\
	&\le c_{13}\|\beta \|^2_{L^\infty(0,T_1;\HH^1)} m_2.
\end{align*}
\item[\mbox{\tiny $\bullet$} \it Step 8:]  {\it The boundedness of the sequence $ ( \int_0^T\|\pg_t u_n\|_2^2\,dt)_n$. }

From equation (\ref{3135}), we  deduce the following inequality for every $n\in\NN$:
\begin{align*}
	\int_0^{T_1}\|\partial_t u_n\|_2^2\,dt& \le
	2\int_0^{T_1}\|Au_n\|_2^2\,dt+\\
&\hh+	2\int_0^{T_1}\|\mathcal{P}_n(- (u _n\!\cdot\!\nabla) u _n -   (\beta\!\cdot\!\nabla) u _n - (u _n\!\cdot\!\nabla)     \beta  +   F   )\|_2^2\,dt\\
&\le  2\int_0^{T_1}\|u_n\|_{\HH^2}^2\,dt+\\
& \hh +8\int_0^{T_1}(\| (u_n\!\cdot\!\nabla )u_n\|_2^2+\| (  \beta\!\cdot\!\nabla) u_n\|_2^2+\| (u_n\!\cdot\!\nabla)    \beta\|_2^2  + \|  F   \|_2^2)\,dt.
\end{align*}
Using this fact along with the boundedness of the sequences in Step 4 to Step 7, we can conclude that the sequence
$(\int_0^{T_1}\|\partial_t u_n\|_2^2\,dt)_n$ is bounded: 	
%x----3150
\begin{align}\label{3150}
\exists m_3>0,\forall n\in\NN,	\int_0^{T_1}\|\partial_t u_n\|_2^2\,dt\le & m_3.
\end{align}	
\item[\mbox{\tiny $\bullet$} \it Step 9:]  {\it  The existence of a convergent subsequence of $(u_n)_n$ and strong solution $u$}.

 Based on the boundedness demonstrated in Step 3 to Step 8, 
we can deduce from compact embeddings that there exists a subsequence of $(u_n)_n$ converging in $L^2(0,T_1,V)$ to a  function  $u$, which is 
 a strong solution for problem (\ref{3007}) and satisfies the properties (\ref{3172})-(\ref{3167}).

The details of the compact embeddings and convergences in this step are presented in Chapter 3 to Chapter 6 in   Robinson \& Rodrigo \& Sadowski \cite{Robinson}.

The proof is complete.\hfill $\Box$
\end{enumerate} 
 
   %%%%%%%%%%%%%%%%%%%%%%%%%%%%%%%%%%%%
\section{The proof of Theorem (\ref{3101})}

%
%\begin{theorem}
%	There exists a bounded domain $\Omega$ in $\RR^3$ of class $C^3$ such that $V$  such that  the problem (\ref{3101}) %with the boundary condition (\ref{3102}) has at least two distinct strong solutions.
%\end{theorem}
%{\bf Proof}.
    In this section, without loss of generality, the radius $r$ is assumed to be $2\pi$. Any other value of $r$ could be used instead, then the function $\sin |x|$ is replaced with $\sin\left(\frac{2\pi}{r}|x|\right)$.
    Additionally, the assumption $u_o\in V\cap \HH^2 \cap C(\bar{\Omega})$ can be replaced with $u_o\in V\cap \HH^2$, since,  due to the Sobolev embedding of $\HH^2(\Omega)$ into $C(\bar{\Omega})$, the function $u_o\in \HH^2$ can be redefined on a set of zero measure in $\Omega$ to  belong to $C(\bar{\Omega})$.
    
     The proof will be divided into several steps.
\begin{enumerate}
\item {\it Recalling the problem $\mbox{(NS)}$ with $r=2\pi$.}

 Let $\Omega$ be the ball  centered at $0$ of the radius $2\pi$ in $\RR^3$.	Suppose $T>0$, $u_o\in V\cap \HH^2 \cap C(\bar{\Omega})$ and $f\in C([0,T];\LL^2)$. 

Consider the Navier-Stokes equations:
%x----------3117
\begin{equation}\label{3117}
(P)\left\{
\begin{array}{rll}
\partial_t u&=\nu \Delta  u- u\!\cdot\!  \nabla  u-\nabla p+ f  &\mbox{ in }\Omega\times [0, T],\\
div\,u&=0&\mbox{ in }\Omega\times [0, T],\\
u(.,0)& = u_o&\mbox{ in }\Omega,\\
\end{array}
\right. 
\end{equation}
with the boundary condition
%x----3144
\begin{equation}\label{3144}
u\cdot\vec{n}=0 \mbox{ on }\pg\Omega.
\end{equation}
\item {\it The first solution $u^*$}.

Applying  Theorem \ref{3110} with $\beta\equiv 0$, we conclude that there exists $T_1\in (0, T]$ such that the problem $(P)$ possesses  a  solution $u^*$  that satisfies
 %x---- 3146,
 \begin{align} 
 	  	\label{3146}
 \begin{array}{ll}
 &\mbox{\tiny $\bullet$ } u|_{t=0}=u_o,\\
 &\dps \mbox{\tiny$\bullet$ } 	u^*\in C([0,T_1];V)\cap L^2(0,T_1;\HH^2), ~\pg_tu\in L^2(0,T_1;\LL^2),\\
 & \dps \mbox{\tiny$\bullet$ }\int_0^s\la 	\partial_t u-\nu \Delta  u+ u\!\cdot\!  \nabla  u - F,\,\varphi\ra\, dt=0  ,
 \end{array}	  	
 \end{align}
 \hh\hh for all test functions $\varphi\in \mathcal{D}_{\sigma,T_1}$ and almost every $s\in [0,T_1]$. 
 
 Based on (\ref{3146}) and the fact $V\subset (H\cap \HH^1)$, we can deduce that
  $u^*$ is a strong solution to the problem $(P)$ on the time interval $[0, T_1]$  in the sense defined in Definition \ref{3127}.

%Furthermore, from the fact $u^*\in C([0,T_1];V)$, it deduces that 
%%x-----3126
%\begin{equation}
%\label{3126}
%\begin{array}{ll}
%&\mbox{\tiny $\bullet$ }\forall t\in [0, T_1], u^*(.,t)\in V.\\
%&\mbox{\tiny $\bullet$ } \forall t\in [0, T_1], u^*(.,t)|_{\pg\Omega}=0.
%\end{array}
%\end{equation}
%
% This solution has the following properties 
%%x-----3126
%\begin{equation}
%\label{3126}
%\left\{\begin{array}{ll}
%& u^*=0 \mbox{ on }\pg\Omega,\\
%&u^*\in C([0,T_1];V)\cap L^2(0,T_1;\HH^2),
%\end{array}
%\right.
%\end{equation}
%and
%%x----3115
%\begin{equation}\label{3115}
% \int_0^s\la \partial_t u^*-\nu \Delta  u^*+ u^*\!\cdot\!  \nabla  u^*- f,\varphi\ra dt=0
%\end{equation}
%for all test functions $\varphi\in \mathcal{D}_\sigma$ and almost every $s\in [0,T_1]$.
%Hence, $u^*$ is a strong solution of the problem $(P_o)$ on the time interval $[0, T_1]$.
\item {\it Introduction of vector fields $w$ and $\beta$.}\\
$-$ Let vector fields $(g_1, g_2, g_3)$, $w$ and $\beta$ be definded by 
%x-----3132
\begin{equation}\label{3132}
\begin{array}{ll}
&\mbox{\tiny $\bullet$ } (g_1(x), g_2(x),g_3(x)):=(\sin |x|,\sin |x|, \sin |x|) \mbox{ for every }x\in \Omega,\\
&\mbox{\tiny $\bullet$ } w:=\nabla\times (g_1, g_2, g_3),\\
&\mbox{\tiny $\bullet$ }\beta(x,t):=u_o(x)+w(x)t\mbox{ for every }(x,t)\in \Omega\times[0,T].
\end{array}
\end{equation}
Here, $x=(x_1,x_2,x_3)$ and $|x|:=\sqrt{x_1^2+x_2^2+x_3^2}$, $g_k(x)\in\RR$ ($k=1,2,3$).\\
$-$
 Taking the curl of the vector field $(g_1, g_2, g_3)$, we obtain
\begin{align*}
\forall x\in\Omega,~& w(x)=\left(\frac{\pg g_3}{\pg x_2}-\frac{\pg g_2}{\pg x_3}\,,\,
	\frac{\pg g_1}{\pg x_3}-\frac{\pg g_3}{\pg x_1}\,,\,
	\frac{\pg g_2}{\pg x_1}-\frac{\pg g_1}{\pg x_2}\right)\\
&=\left( \cos|x|\frac{x_2}{|x|}-\cos|x|\frac{x_3}{|x|}\,,\,
\cos|x|\frac{x_3}{|x|}-\cos|x|\frac{x_1}{|x|}\,,\,
\cos|x|\frac{x_1}{|x|}-\cos|x|\frac{x_2}{|x|}
\right)	
\end{align*}

For every $x\in\pg\Omega$, we have $\cos|x|=cos (2\pi)=1$, and therefore
\begin{align*}
w(x)\cdot \vec{n}(x)&=\left(\frac{x_2-x_3}{r}, \frac{x_3-x_1}{r}, \frac{x_1-x_2}{r}\right)
\cdot \left(\frac{x_1}{r},\frac{x_2}{r}, \frac{x_3}{r},\right)\\
&=\frac{x_1(x_2-x_3)+x_2(x_3-x_1)+x_3(x_1-x_2)}{r^2}=0.
\end{align*} 	

Hence, the properties of $w$ on the boundary $\pg\Omega$ are as follows:
%x----3133
\begin{align}\label{3133}
	\left\{
		\begin{array}{ll}
& w\not\equiv 0\mbox{ on }\pg\Omega \\
& w\cdot \vec{n}=0  \mbox{ on }\pg\Omega.
	\end{array}
	\right.
\end{align}
$-$
Since the divergence of a curve is zero, we have $div(\nabla\times (g_1, g_2, g_3))=0$,  and therefore
%x-----3131
\begin{equation}
\label{3131}
div~w=0 \mbox{ in }\Omega.
\end{equation} 
%$-$ Denote $h(t):=w t$. From (\ref{3132}), (\ref{3133}), (\ref{3131}), one obtains
%\begin{align}
%h \in C([0, T], (H\cap \HH^2)\backslash V). 
%\end{align}
$-$ From (\ref{3132})-(\ref{3131}), we obtain
%x-----3124
\begin{equation}
\label{3124}
\left\{
\begin{array}{ll}
w\in (H\cap \HH^2)\backslash V\\
\beta\in C^1([0,T]; H\cap \HH^2)\\
\pg_t\beta =w.
\end{array}
\right.
\end{equation}
$-$ Furthermore, from (\ref{3124}) we deduce that
\begin{align*}
	 \int_0^{T}\!\!\ino |\beta|^2|\nabla \beta|^2\,dx\,dt&\le \int_0^{T}\!\!\ino |\nabla \beta|^2\,dx\|\beta(.,t)\|_{\LL^\infty}^2\,dt\\
	 &\le \| \beta\|_{L^\infty(0,T;\HH^1)}^2\int_0^T\|\beta(.,t)\|_{\LL^\infty}^2\,dt\\
	 &	\le c\| \beta\|_{L^\infty(0,T;\HH^1)}^2\int_0^T\|\beta\|_{\HH^2}^2\,dt<\infty.
	% \\ 
	 % &\le c\| \beta\|_{L^\infty(0,T;\HH^1)}^2\| \beta\|_{L^\infty(0,T;\HH^2)}^2<+\infty
\end{align*}	
 
 Therefore,  
\begin{align}\label{3164}
	 (\beta\cdot\nabla) \beta \in L^2(0,T;\LL^2).
\end{align}	
\item {\it The second solution $\tilde{u}$}.\\
$-$
Denote 
%x-----3119
\begin{equation}\label{3119}
F:=	-\partial_t  \beta+ \nu\Delta \beta 	     -    (\beta \cdot\!\nabla) \beta +f.
\end{equation}

Based on this fact, (\ref{3124}), (\ref{3164}) and the assumption 
$f\in  L^2(0, T; \LL^2)$, we deduce that $F \in L^2(0, T; \LL^2)$.

Consider the problem for unknown $(v,p)$:
  	%x----------3116
  	\begin{equation}\label{3116}
  	\left\{
  	\begin{array}{rll}
  	\partial_t v&=\nu \Delta  v- v\!\cdot\!  \nabla  v- \beta \!\cdot\!  \nabla  v- v\!\cdot\!  \nabla  \beta -\nabla p+ F&\mbox{ in }\Omega,\\
  	div~v&=0&\mbox{ in }\Omega,\\
  	v(.,0)& = 0&\mbox{ in }\Omega.\\
  	\end{array}
  	\right. 
  	\end{equation}
   
  By Theorem \ref{3110}, there exist $T_2\in (0,T_1]$ such that   the problem $(\ref{3116})$ possesses a  solution $v$ that satisfies
  %x-----3149, 3170, 3120
  \begin{align} \label{3149} 
     &  \mbox{\tiny$\bullet$ }  v|_{t=0}=0,\\ 	
 \label{3170}
  	&  \mbox{\tiny$\bullet$ }\dps v\in C([0,T_2];V)\cap L^2(0,T_2;\HH^2), \pg_tv\in L^2(0,T_2;\LL^2),\\
  \label{3120}
  &	\mbox{\tiny$\bullet$ } \int_0^s\la  \partial_t v-\nu \Delta  v+ v\!\cdot\!  \nabla  v+ \beta \!\cdot\!  \nabla  v+ v\!\cdot\!  \nabla  \beta - F,~ \varphi\ra dt=0 
   \end{align}
  for all test functions $\varphi\in \mathcal{D}_{\sigma, T_2}$ and almost every $s\in [0,T_2]$.
 \\ 
  $-$ Let $\tilde{u}$ be the vector field  defined by
 %x------3109
 \begin{equation}\label{3109}
 \tilde{u}:=v+ \beta \mbox{ for every } (x,t)\in\Omega\times[0, T_2]. 
  \end{equation}
  
   Based on (\ref{3109}), the properties of $v$ (as stated in  (\ref{3170})),  and\\
   the properties of $\beta$ (as stated in  (\ref{3124})), we deduce that 
   %x-----3122
   \begin{equation}\label{3122}
   \tilde{u}\in  C([0,T_2];H\cap \HH^1)\cap L^2(0,T_2;\HH^2)\mbox{ and }\pg_t\tilde{u}\in L^2(0,T_2;\LL^2).
     \end{equation} 
%  From this fact and Definition \ref{3127} of the space $H$, we deduce
 %x----3163
 %\begin{equation}
 %\label{3163}
%div~\tilde{u}=0 \mbox{ in }\Omega, \mbox{ and } \tilde{u}\cdot \vec{n}=0 \mbox{ on }\pg\Omega.
 %\end{equation}   

 Based on (\ref{3109}), (\ref{3132}),  and (\ref{3149}),  we deduce that: 
   %x----3011
   \begin{equation}\label{3011}
  \tilde{u}|_{t=0}=v|_{t=0}+u_o+[w.t]_{t=0}=0+u_o+0=u_o.
   \end{equation}
   
  By (\ref{3109}), we have $v=	\tilde{u}- \beta$. 			
  Substituting $v$ with $\tilde{u}- \beta $ in  equation (\ref{3120}), and using (\ref{3119}), we obtain the following equation for all test functions    $\varphi\in \mathcal{D}_{\sigma,T_2}$ and almost every $s\in [0,T_2]$:
  \begin{align*}
 % 	\int_0^s\la  \partial_t v-\nu \Delta  v+ v\!\cdot\!  \nabla  v+ \beta \!\cdot\!  \nabla  v+ v\!\cdot\!  \nabla  \beta - F, \varphi\ra dt&=0 \\
    \int_0^s\la 	\partial_t(\tilde{u}-\!\! \beta )-\nu \Delta (\tilde{u}- \!\!\beta )+(\tilde{u}- \!\!\beta )\!\cdot\!  \nabla (\tilde{u}- \!\!\beta )+\beta \cdot \nabla(\tilde{u}-\beta),\varphi\ra dt+\hh&\\
    +\int_0^s \la(\tilde{u}-\beta) \cdot \nabla \beta- F,\varphi\ra\, dt &=0\\
   	\int_0^s \la 	\displaystyle \partial_t \tilde{u} - \nu\Delta \tilde{u}+ \tilde{u}\!\cdot\!\nabla \tilde{u}  
  +\underbrace{-	\partial_t  \beta+\nu\Delta \beta 	- \beta \!\cdot\!\!\nabla \beta 
  }_{F-f}-F, \varphi \ra\, dt&=0\\
  \int_0^s \la 	\displaystyle \partial_t \tilde{u} - \nu\Delta \tilde{u}+ \tilde{u}\!\cdot\!\nabla \tilde{u}  
  -f, \varphi \ra \,dt&=0.
    \end{align*}	
 
 Based on this fact, (\ref{3011}) and (\ref{3122}), we deduce that  $\tilde{u}$ is a strong solution to the problem $(P)$ on the time interval $[0, T_2]$
  in the sense defined in Definition  \ref{3127}.
\item {\it Two distinct solutions. }

Let us prove that two solution $u^*$ and $\tilde{u}$ are distinct.

Due to  (\ref{3109}) and (\ref{3132}), we have $\tilde{u}=v+u_o+wt$. From this fact,
the property $v(.,t)\in V$ (due to (\ref{3170})), the assumption $u_o\in V$, and $w|_{\pg\Omega}\not\equiv 0$ (as stated in (\ref{3133})),  it implies that
  \begin{align*}
 \forall t\in (0, T_2],	\tilde{u}(.,t)|_{\pg\Omega}=v(.,t)|_{\pg\Omega}+u_o|_{\pg\Omega}+wt|_{\pg\Omega}
 =0+0+wt|_{\pg\Omega}\not\equiv 0,
 	\end{align*}
 and therefore  $\tilde{u}(.,t)\not\in V$.

 Besides, by (\ref{3146}), $u^*(.,t)\in V$  for every $t\in [0, T_2]$.

 Hence, $u^*$ and $\tilde{u}$ are two distinct strong solutions to the problem $(P)$ on the time interval $[0, T_2]$. 
\end{enumerate}	

The proof is complete.\hfill 
$\Box$	
            %%%%%%%%%%%%%%%%%%%%%%
% % % % % % % % % % % % % % % % % % % % % % % % % % % %
% % % % % % % % % % % % % % % % %========================
%%%%%%%%%%%%%%%

\end{document}